\documentclass[preprint]{elsarticle}
\usepackage[letterpaper,top=2cm,bottom=2cm,left=3cm,right=3cm,marginparwidth=1.75cm]{geometry}
\usepackage{CJKutf8}
\usepackage{amsmath, amsfonts, amssymb} 
\usepackage[english]{babel}
\allowdisplaybreaks[4]
\usepackage{lipsum}
\usepackage{graphicx}   
\usepackage{url}        
\usepackage{bm}         
\usepackage{fancyhdr} 
\usepackage{multirow}
\usepackage{booktabs}
\usepackage{algorithm}
\usepackage{algorithmic}
\usepackage{setspace}
\usepackage{multicol}
\usepackage{graphicx}
\usepackage[colorlinks=true, allcolors=blue]{hyperref}
\journal{arXiv}
\begin{document}
\begin{frontmatter}
\title{Uniqueness of mild solutions to the Navier-Stokes equations in weak-type $L^d$ space 
}
\author[first]{Zhirun Zhan}
\affiliation[first]{organization={Division of Mathematics and Mathematical Sciences, Kyoto University},
            addressline={Oiwakecho, Kitashirakawa, Sakyo-ku}, 
            city={Kyoto},
            postcode={606-8502}, 
            country={Japan}}
\ead{zhan.zhirun.82z@st.kyoto-u.ac.jp}
\date{}
\begin{abstract}
This paper deals with the uniqueness of mild solutions to the forced or unforced Navier-Stokes equations in the whole space. It is known that the uniqueness of mild solutions to the unforced Navier-Stokes equations holds in $L^{\infty}(0,T;L^d(\mathbb{R}^d))$ when $d\geq 4$, and in $C([0,T];L^d(\mathbb{R}^d))$ when $d\geq3$. As for the forced Navier-Stokes equations, when $d\geq3$ the uniqueness of mild solutions in $C([0,T];L^{d,\infty}(\mathbb{R}^d))$ with force $f$ and initial data $u_{0}$ in some proper Lorentz spaces is known. In this paper we show that for $d\geq3$, the uniqueness of mild solutions to the forced Navier-Stokes equations in $ C((0,T];\widetilde{L}^{d,\infty}(\mathbb{R}^d))\cap L^\beta(0,T;\widetilde{L}^{d,\infty}(\mathbb{R}^d))$ for $\beta>2d/(d-2)$ holds when there is a mild solution in $C([0,T];\widetilde{L}^{d,\infty}(\mathbb{R}^d))$ with the same initial data and force. Here $\widetilde{L}^{d,\infty}$ is the closure of ${L^{\infty}\cap L^{d,\infty}}$ with respect to $L^{d,\infty}$ norm.
\end{abstract}
\end{frontmatter}
Declarations of interest: none.
\section{Introduction}
In this paper we consider the following Navier-Stokes equations in $(0,T)\times\mathbb{R}^d$ for $d\geq3$:
\begin{equation}
\begin{cases}
    \partial_{t}u=\Delta u-(u\cdot \nabla)u- \nabla p+f \quad $in    $(0,T)\times\mathbb{R}^d,\\
 \nabla\cdot u=0 \quad    $in    $[0,T)\times\mathbb{R}^d,\\
 u(0,x)=u_{0}(x) \quad   $in    $\mathbb{R}^d.\\
\end{cases}
\end{equation}
Here $T\in(0,\infty)$, $u=(u_{1}(t,x),...,u_{d}(t,x))$ is the unknown velocity field, $p=p(t,x)$ is the unknown pressure field, $f=(f_{1}(t,x),...,f_{d}(t,x))$ is a given force, and $u_{0}$ is a given initial velocity field. We define 
$\nabla p=(\partial_{x_1}p,...,\partial_{x_d}p)$ and $\nabla\cdot u=\sum_{j=1}^{d}\partial_{x_j}u_{j}$. To shorten the notation, we omit the superscript of vector space in this paper, for example, we note $(L^p(\mathbb{R}^d))^d$ as $L^p(\mathbb{R}^d)$. When $f=0$, the uniqueness of mild solutions to the unforced Navier-Stokes equations in $C([0,T];L^3(\mathbb{R}^3))$ has been proved by Furioli, Lemari$\acute{e}$-Rieusset, and Terraneo \cite{furioli2000unicite}. In \cite{furioli2000unicite} the key step is to show the bilinear estimate for the fluctuation $u-e^{t\Delta}u_{0}$ in the homogeneous Besov space $\dot{B}^{1/2,\infty}_{2}(\mathbb{R}^3)$. Meyer \cite{meyer1996wavelets} gave the alternative proof of the uniqueness in $C([0,T];L^3(\mathbb{R}^3))$ by establishing the bilinear estimate in the Lorentz space $L^{3,\infty}(\mathbb{R}^3)$. Monniaux \cite{monniaux1999uniqueness} gave a simpler proof based on the maximal $L^pL^q$ regularity property for the heat kernel. Gala \cite{gala2007note} also gave a proof based on the bilinear estimate in $L^4(0,T;L^3(\mathbb{R}^3))$ and the Grönwall's inequality. The uniqueness in $L^{\infty}(0,T;L^3(\mathbb{R}^3))$ remains an open problem, but one can show that the uniqueness of mild solutions to the unforced Navier-Stokes equations in $L^{\infty}(0,T;L^2(\mathbb{R}^3))\cap L^{2}(0,T;\dot{H}^1(\mathbb{R}^3))\cap L^{\infty}(0,T;L^3(\mathbb{R}^3))$ holds with initial data $u_{0}\in L^2(\mathbb{R}^3)\cap L^3(\mathbb{R}^3)$, following the energy method given by Von Wahl \cite{von1985equations}. Moreover, the uniqueness of mild solutions to the unforced Navier-Stokes equations in $L^{\infty}(0,T;L^d(\mathbb{R}^d))$ for $d\geq4$ was proved by Lions and Masmoudi \cite{lions1998unicite}. The uniqueness of mild solutions to the unforced Navier-Stokes equations in $C([0,T];bmo^{-1})\cap L_{loc}^{\infty}((0,T);L^{\infty}(\mathbb{R}^d))$ with initial data in $vmo^{-1}$ was proved by Miura \cite{miura2005remark}. Besides, when $d=2$, Brezis \cite{ben1994global} proved the uniqueness when the vorticity field belongs to  $C([0,\infty);L^1(\mathbb{R}^2))\cap C((0,\infty);L^{\infty}(\mathbb{R}^2))$. As for the forced Navier-Stokes equations, Farwig,  Nakatsuka, and Taniuchi \cite{farwig2015uniqueness} proved the uniqueness of mild solutions in $C([0,T];\widetilde{L}^{3,\infty}(\mathbb{R}^3))$. Here $\widetilde{L}^{d,\infty}$ is the closure of ${L^{\infty}\cap L^{d,\infty}}$ with respect to $L^{d,\infty}$ norm. Yamazaki \cite{yamazaki2000navier} proved the uniqueness of mild solutions in $BC((0,\infty);L^{d,\infty}(\mathbb{R}^d))$ for $f=\nabla\cdot F, F\in BC((0,\infty);L^{d/2,\infty})$ with the smallness assumption of solutions and force in this space when $d\geq3$. Here ``$BC$'' stands for ``bounded and continuous''. Okabe and Tsutsui \cite{okabe2021remark} proved the uniqueness of mild solutions in $C([0,T];L^{d,\infty}(\mathbb{R}^d))$ with $u_{0}\in \widetilde{L}^{d,\infty}(\mathbb{R}^d)$, $f\in C([0,T];L^1(\mathbb{R}^3))$ for $d=3$ and $f\in C([0,T];\widetilde{L}^{d/3,\infty}(\mathbb{R}^d))$ for $d\geq4$.

Apart from above, there have been many results about the unconditional uniqueness in $C([0,T];X)$ for some appropriate Banach space $X$ (see \cite{2011The} for the Herz spaces, and \cite{ferreira2016bilinear,ferreira2020bilinear} for the Morrey-Lorentz spaces). However, the uniqueness in $L^{\infty}(0,T;X)$ remains an open problem in general. In fact, the situation changes drastically when the solution orbit $\{u(t)\}_{t\in[0,T]}$ is not assumed to be compact nor small in scale critical spaces. Indeed, recently, Albritton, Bru$\acute{e}$, and Colombo \cite{albritton2022non} have proved the non-uniqueness of the Leray-Hopf solutions for $d=3$. More precisely, in \cite{albritton2022non}, two distinct weak solutions $u,v\in L^{\infty}(0,T;L^{2}(\mathbb{R}^3))\cap L^{2}(0,T;\dot{H}^{1}(\mathbb{R}^3))$ to the Navier-Stokes equations are constructed for some $T\in(0,\infty)$, initial data $u_{0}=0$, and $f\in L^1(0,T;L^2(\mathbb{R}^3))$. We can check that the solutions constructed in \cite{albritton2022non} belong to $BC((0,T];L^3(\mathbb{R}^3))$. This implies that the continuity up to $t=0$ of the solution, or the compactness of the solution orbit $\{u(t)\}_{t\in[0,T]}$, in the scale critical space such as $L^3(\mathbb{R}^3)$ is crucial for determining the uniqueness/non-uniqueness of mild solutions to the Navier-Stokes equations. Summarizing the known results, we see that achieving the uniqueness theorem becomes highly nontrivial in the following cases:\\
(\romannumeral1) When the solution orbit $\{u(t)\}_{t\in[0,T]}$ is not assumed to be compact nor small in scale critical spaces.\\
(\romannumeral2) When the space dimension is low, i.e., d=2,3. For example, when $f=0$, the unconditional uniqueness in $L^{\infty}(0,T;L^d(\mathbb{R}^d))$ is still open when $d=2,3$, though $d\geq4$ is settled. In fact, the unconditional uniqueness in $C([0,T];L^d(\mathbb{R}^d))$ is open when $d=2$, though $d\geq3$ is settled. \\
(\romannumeral3) When the nontrivial force is taken into account. For example, the uniqueness in $BC((0,T];L^3(\mathbb{R}^3))$ breaks down with some $f$ large in a scale critical space. 

 The aim of this paper is to contribut to the uniqueness theory in such delicate cases for $d\geq3$. Indeed, in this paper we show that, without any smallness condition, the uniqueness of mild solutions to the Navier-Stokes equations holds in $C((0,T];\widetilde{L}^{d,\infty}(\mathbb{R}^d))\cap L^\beta(0,T;\widetilde{L}^{d,\infty}(\mathbb{R}^d))$ for $\beta>2d/(d-2)$, if there is a mild solution in $C([0,T_{0}];\widetilde{L}^{d,\infty}(\mathbb{R}^d))$ for some $T_{0}\in(0,T]$ with the same initial data and force. We note that when $u\in C((0,T];\widetilde{L}^{d,\infty}(\mathbb{R}^d))\cap L^\beta(0,T;\widetilde{L}^{d,\infty}(\mathbb{R}^d))$ for $\beta>2d/(d-2)$, the solution orbit $\{u(t)\}_{t\in[0,T]}$ is not necessarily compact (even in the case $\beta=\infty$) in $\widetilde{L}^{d,\infty}$.  In particular, our result significantly extends the uniqueness results of \cite{furioli2000unicite} for the case $f=0$ and \cite{okabe2021remark} for the case $f\neq0$, since the existence of the mild solution in $C([0,T_0];\widetilde{L}^{d,\infty})$ is known by \cite{okabe2021remark} when $u_0\in X_{\sigma}^{d,\infty}$, $f\in C([0,T];L^1(\mathbb{R}^3))$ for $d=3$ and $f\in C([0,T];\widetilde{L}^{d/3,\infty}(\mathbb{R}^d))$ for $d\geq4$. Here $X_{\sigma}^{d,\infty}$ is the maximal subspace of $L_{\sigma}^{d,\infty}(\mathbb{R}^d)$ such that the Stokes semigroup is a $C_0$-semigroup. Note that the existence of the mild solution in $C([0,T];L^d(\mathbb{R}^d))$ for $u_0\in L^d(\mathbb{R}^d)$ and $f=0$ is due to \cite{kato1984strong}. It should be emphasized here that, when $d=3$, our results gives the first unconditional uniqueness without assuming the compactness nor the smallness for the solution orbit $\{u(t)\}_{t\in[0,T]}$ in scale critical spaces, even for the simplest case $u_0\in L^3(\mathbb{R}^3)$ and $f=0$. Moreover, as explained later, our approach is rather generic and robust in nature. Hence, it will be useful to investigate the unconditional uniqueness of various nonlinear PEDs in a scale critical space without smallness nor compactness of the solution orbit.
To state our result precisely, we first introduce the following definitions.\\[5pt]
\textbf{Definition 1} (Lorentz spaces)\textbf{.} It is well known that the Lorentz space can be realized as a real interpolation space $L^{p,q}=(L^{p_{0}},L^{p_{1}})_{\theta,q}$ with $1/p=(1-\theta)/p_{0}+\theta/p_{1}$, $1\leq p_{0},p_{1},q\leq\infty$ and $\theta\in(0,1)$. We define the divergence free Lorentz space $L_{\sigma}^{p,q}:=\mathbb{P}L^{p,q}$ with $p\in(1,\infty)$, $q\in[1,\infty]$. Here $\mathbb{P}$ is the Leray projection $\mathbb{P}=I-\nabla\Delta^{-1}\nabla\cdot$.\\[5pt]
\textbf{Remark 1.} We recall the Helmholtz decomposition for $L^p$ in the whole space $\mathbb{R}^d$: 
$$L^p(\mathbb{R}^d)=L^p_{\sigma}(\mathbb{R}^d)\oplus G^p(\mathbb{R}^d), \quad p\in(1,\infty),$$ where $L^p_{\sigma}(\mathbb{R}^d)=\{f\in L^p(\mathbb{R}^d)| \nabla\cdot f=0\}=\mathbb{P}L^p$. Since $\mathbb{P}$ is surjective and bounded from $L^{p}$ to $L^{p}_{\sigma}$, we notice that $\mathbb{P}L^{p,q}(\mathbb{R}^d)=(\mathbb{P}L^{p_{0}}(\mathbb{R}^d),\mathbb{P}L^{p_{1}}(\mathbb{R}^d))_{\theta,q}$ with $1/p=(1-\theta)/p_{0}+\theta/p_{1}$, $p_{0},p_{1}\in(1,\infty)$ and $q\in[1,\infty]$ (see \cite{okabe2021remark}). Thus, $L^{p,q}_{\sigma}(\mathbb{R}^d)=(L_{\sigma}^{p_{0}}(\mathbb{R}^d),L_{\sigma}^{p_{1}}(\mathbb{R}^d))_{\theta,q}$. Moreover, we obtain the Helmholtz decomposition for $L^{p,q}$ by the real interpolation, 
$$L^{p,q}(\mathbb{R}^d)=L^{p,q}_{\sigma}(\mathbb{R}^d)\oplus G^{p,q}(\mathbb{R}^d), \quad p\in(1,\infty),q\in[1,\infty],$$
where $G^{p,q}(\mathbb{R}^d)=\{u=\nabla p\in L^{p,q}(\mathbb{R}^d)|p\in L^p_{loc}\}$ (see \cite{yamazaki2000navier}).\\[5pt] 
\textbf{Definition 2.} We also introduce the following spaces for $p\in(1,\infty)$. $$\widetilde{L}^{p,\infty}(\mathbb{R}^d)=\overline{L^{p,\infty}(\mathbb{R}^d)\cap L^{\infty}(\mathbb{R}^d)}^{L^{p,\infty}}, \widetilde{L}_{\sigma}^{p,\infty}=\overline{L_{\sigma}^{p,\infty}(\mathbb{R}^d)\cap L^{\infty}(\mathbb{R}^d)}^{L^{p,\infty}}.$$
\textbf{Remark 2.} $\widetilde{L}^{p,\infty}(\mathbb{R}^d)=\overline{L^{p,\infty}(\mathbb{R}^d)\cap L^{\infty}(\mathbb{R}^d)}^{L^{p,\infty}}=\overline{L^{p,\infty}(\mathbb{R}^d)\cap L^{s}(\mathbb{R}^d)}^{L^{p,\infty}}$ for any $s\in(d,\infty)$. Indeed, $L^{p,\infty}(\mathbb{R}^d)\cap L^{\infty}(\mathbb{R}^d)=L^{p,\infty}(\mathbb{R}^d)\cap L^{\infty}(\mathbb{R}^d)\cap L^{s}(\mathbb{R}^d)\subset L^{p,\infty}(\mathbb{R}^d)\cap L^{s}(\mathbb{R}^d)$. By the $K$-method of real interpolation, any $u\in L^{s}$ can be decomposed to $u_{1}+u_{2}$ with $u_{1}$ small in $L^{d,\infty}$ and $u_{2}\in L^{\infty}$. Therefore, $\widetilde{L}^{p,\infty}(\mathbb{R}^d)=\overline{L^{p,\infty}(\mathbb{R}^d)\cap L^{s}(\mathbb{R}^d)}^{L^{p,\infty}}$ for any $s\in(d,\infty]$.\vspace{5pt}

We also want to introduce some Sobolev spaces equipped with the Lorentz norm.\\[5pt]
\textbf{Definition 3} (Second order Sobolev-Lorentz spaces)\textbf{.}
$$W^{2,p,q}(\mathbb{R}^d)=\{u|u,\nabla u,\nabla^{2}u\in L^{p,q}(\mathbb{R}^d)\},\quad p\in(1,\infty),q\in[1,\infty],$$
$$W^{2,p,q}_{\sigma}(\mathbb{R}^d)=\mathbb{P}W^{2,p,q}(\mathbb{R}^d),\quad p\in(1,\infty),q\in[1,\infty].$$
\textbf{Remark 3.} (\romannumeral1) Naturally, we want to figure out if this space is equivalent to the interpolation space $\widehat{W}^{2,p,q}(\mathbb{R}^d)=(W^{2,p_{0}}(\mathbb{R}^d),W^{2,p_{1}}(\mathbb{R}^d))_{\theta,q}$ with $1/p=(1-\theta)/p_{0}+\theta/p_{1},\theta\in(0,1),p_{0},p_{1}\in(1,\infty),q\in[1,\infty].$ Actually, these two definitions are equivalent. By the $J$-method of real interpolation, we easily check that $\widehat{W}^{2,p,q}\subset W^{2,p,q}$. As a classic result, the resolvent $(\lambda I-\Delta)^{-1}$ of the Laplace operator is bijective and bounded from $L^{p}(\mathbb{R}^d)$ to $W^{2,p}(\mathbb{R}^d)$ with $\lambda$ positive and $p\in(1,\infty)$. By the $J$-method again, we obtain that $(\lambda I-\Delta)^{-1}$ is bijective and bounded from $L^{p,q}(\mathbb{R}^d)$ to $\widehat{W}^{2,p,q}(\mathbb{R}^d)$. Assume that $u\in W^{2,p,q}(\mathbb{R}^d)$, set $f=(\lambda-\Delta)u$, we have $v=(\lambda I-\Delta)^{-1}f\in\widehat{W}^{2,p,q}(\mathbb{R}^d)$ and $(\lambda-\Delta)(u-v)=0$ in $L^{p,q}(\mathbb{R}^d)$. Since $(\lambda I-\Delta)^{-1}$ is bounded and bijective on the Schwartz space with $\lambda$ positive, we derive $u-v=0$ in $S'(\mathbb{R}^d)$, which implies $u=v$ in $\widehat{W}^{2,p,q}(\mathbb{R}^d)$. Therefore, $W^{2,p,q}(\mathbb{R}^d)=\widehat{W}^{2,p,q}(\mathbb{R}^d)$.\\[5pt]
(\romannumeral2) Since $\mathbb{P}$ is bounded and surjective from $W^{2,p,q}(\mathbb{R}^d)$ to $W^{2,p,q}_{\sigma}(\mathbb{R}^d)$, we have $$W^{2,p,q}_{\sigma}(\mathbb{R}^d)=\mathbb{P}W^{2,p,q}(\mathbb{R}^d)=(\mathbb{P}W^{2,p_{0}}(\mathbb{R}^d),\mathbb{P}W^{2,p_{1}}(\mathbb{R}^d))_{\theta,q}=(W_{\sigma}^{2,p_{0}}(\mathbb{R}^d),W_{\sigma}^{2,p_{1}}(\mathbb{R}^d))_{\theta,q},$$
where $1/p=(1-\theta)/p_{0}+\theta/p_{1},\theta\in(0,1),p_{0},p_{1}\in(1,\infty),q\in[1,\infty],$ and $W_{\sigma}^{2,p}(\mathbb{R}^d)=\mathbb{P}W^{2,p}(\mathbb{R}^d)=\{u\in W^{2,p}(\mathbb{R}^d)|\nabla\cdot u=0\}$. We notice that $\mathbb{P}\Delta u=\Delta \mathbb{P}u$.\\[5pt]
(\romannumeral3) The space $C_{c,\sigma}^{\infty}(\mathbb{R}^d)=\{u\in C_{c}^{\infty}(\mathbb{R}^d),\nabla\cdot u=0\}$ is dense in $W_{\sigma}^{2,p}(\mathbb{R}^d)$ for $p\in(1,\infty)$. Therefore, we obtain that $C_{c,\sigma}^{\infty}(\mathbb{R}^d)$ is dense in $W_{\sigma}^{2,p,q}(\mathbb{R}^d)$ for $p\in(1,\infty)$, $q\in[1,\infty)$ by the real interpolation method. Moreover, $C_{c,\sigma}^{\infty}((0,T)\times\mathbb{R}^d)=\{u\in C_{c}^{\infty}((0,T)\times\mathbb{R}^d),\nabla_x\cdot u =0\}$ is dense in $L^{p^*}(0,T;W_{\sigma}^{2,p,q}(\mathbb{R}^d))$ for $p\in(1,\infty)$ and $p^*,q\in[1,\infty)$, since $L^{p^*}(0,T;W_{\sigma}^{2,p,q}(\mathbb{R}^d))$ is a Bochner space.\\[5pt]
(\romannumeral4) The reason we introduce the Sobolev spaces with Lorentz norm is that our results are stated on the scale critical Lorentz space $L^{d,\infty}(\mathbb{R}^d)$. We notice that the domain of the Stokes operator in $L^{d,\infty}(\mathbb{R}^d)$ is $D(\mathbb{P}\Delta)=W_{\sigma}^{2,d,\infty}(\mathbb{R}^d)$.
According to the classical abstract semigroup theory stated in \cite{lunardi2012analytic}, for $f\in L_{\sigma}^{d,\infty}(\mathbb{R}^d)$, $\lim_{t\to0^+}\|e^{t\Delta}f-f\|_{L^{d,\infty}}=0$ if and only if $f\in \overline{D(\mathbb{P}\Delta)}^{L^{d,\infty}}$. Set $X_{\sigma}^{d,\infty}=\overline{D(\mathbb{P}\Delta)}^{L^{d,\infty}}$ as the maximal subspace of $L_{\sigma}^{d,\infty}(\mathbb{R}^d)$ such that the Stokes semigroup is a $C_{0}$-semigroup on $X_{\sigma}^{d,\infty}$. We have the following conclusion (see \cite{okabe2021remark}), 
$$\overline{C_{c,\sigma}^{\infty}(\mathbb{R}^d)}^{L^{d,\infty}}\subsetneqq  X_{\sigma}^{d,\infty}\subsetneqq \widetilde{L}^{d,\infty}_{\sigma}(\mathbb{R}^d)\subsetneqq L^{d,\infty}_{\sigma}(\mathbb{R}^d).$$
\textbf{Definition 4} (mild solution)\textbf{.} Let $d\geq3$, $\beta\geq2$, $u_{0}\in \widetilde{L}^{d,\infty}_{\sigma}(\mathbb{R}^d)$, and $f\in L^1(0,T;L^{1}(\mathbb{R}^d))$ or $f\in L^1(0,T;L^{r,\infty}(\mathbb{R}^d))$ for some $r\in (1,\infty)$. For $u\in C((0,T];\widetilde{L}_{\sigma}^{d,\infty}(\mathbb{R}^d))\cap L^\beta(0,T;\widetilde{L}_{\sigma}^{d,\infty}(\mathbb{R}^d))$, we call $u$ is a mild solution of $(1)$ if $u$ satisfies
\begin{equation}
  u(t)=e^{t\Delta}u_{0}+\int_{0}^{t}e^{(t-s)\Delta}\mathbb{P}(f-\nabla\cdot(u\otimes u))ds,\quad t\in(0,T).
\end{equation}
Here $e^{t\Delta}$ is the heat semigroup in the whole space. \vspace{5pt}

We give another definition of solutions. Readers can confirm that these two definitions are actually equivalent by Proposition 1 in Section 2.\\[5pt]
\textbf{Definition 5} (very weak solution)\textbf{.} Let $d\geq3$, $\beta\geq2$, $u_{0}\in \widetilde{L}^{d,\infty}_{\sigma}(\mathbb{R}^d)$, and $f\in L^1(0,T;L^{1}(\mathbb{R}^d))$ or $f\in L^1(0,T;L^{r,\infty}(\mathbb{R}^d))$ for some $r\in (1,\infty)$. Let $u\in C((0,T];\widetilde{L}_{\sigma}^{d,\infty}(\mathbb{R}^d))\cap L^\beta(0,T;\widetilde{L}_{\sigma}^{d,\infty}(\mathbb{R}^d))$. We call such $u$ is a very weak solution of $(1)$ if the following equality holds for any $\phi\in (C_{c,\sigma}^{\infty}([0,T]\times \mathbb{R}^d))^d$ satisfying $\phi(T)\equiv0$.
\begin{equation}
\int_{\mathbb{R}^d}u_{0}\cdot\phi(0,x)dx+\int_{0}^{T}\int_{\mathbb{R}^d}f\cdot\phi+u\cdot(\partial_{t}\phi+\Delta\phi)+\sum_{1\leq i,j\leq d}u_{i}u_{j}\partial_{i}\phi_{j}dxdt=0.   
\end{equation}
\vspace{5pt}

Our main result is stated as follows:\\[5pt]
\textbf{Theorem 1.} Let $d\geq3$, $\beta>2d/(d-2)$, $u_{0}\in \widetilde{L}^{d,\infty}_{\sigma}(\mathbb{R}^d)$, and $f\in L^1(0,T;L^{1}(\mathbb{R}^d))$ or $f\in L^1(0,T;L^{r,\infty}(\mathbb{R}^d))$ for some $r\in (1,\infty)$. Assume that $u_{1},u_{2}\in C((0,T];\widetilde{L}_{\sigma}^{d,\infty}(\mathbb{R}^d))\cap L^\beta(0,T;\widetilde{L}_{\sigma}^{d,\infty}(\mathbb{R}^d))$ are two mild solutions of $(1)$ in $(0,T)\times\mathbb{R}^d$ with the same force $f$ and initial data $u_{0}$. If there exists a mild solution $v\in C([0,T_{0}];\widetilde{L}^{d,\infty}_{\sigma}(\mathbb{R}^d))$ of (1) in $(0,T_{0})\times\mathbb{R}^d$ for some $T_{0}\in(0,T]$, with the same force $f$ and initial data $u_{0}$, then $u_{1}=u_{2}$ in $(0,T)\times\mathbb{R}^d$.\vspace{5pt}

It is known in \cite{okabe2021remark} that for initial data $u_{0}\in X_{\sigma}^{d,\infty}$ and $f\in C([0,T];L^1(\mathbb{R}^3))$ when $d=3$, or $f\in C([0,T];\widetilde{L}^{d/3,\infty}(\mathbb{R}^d))$ when $d\geq4$, there exists a mild solution $v\in C([0,T_{0}];\widetilde{L}_{\sigma}^{d,\infty}(\mathbb{R}^d))$ of $(1)$ in $(0,T_{0})\times\mathbb{R}^d$ for some $T_{0}\in(0,T]$. We can also construct such mild solution $v$ with $f$, which satisfies $t^{\alpha}f\in L^\infty(0,T;L^{d,\infty}(\mathbb{R}^d))$ for some $\alpha<1$. Therefore, we obtain the following corollary.\\[5pt]
\textbf{Corollary 1.} (\romannumeral1) Let $d\geq3$, $\beta>2d/(d-2)$, $u_{0}\in X_{\sigma}^{d,\infty}$, $f\in C([0,T];L^1(\mathbb{R}^3))$ for $d=3$, or $f\in C([0,T];\widetilde{L}^{d/3,\infty}(\mathbb{R}^d))$ for $d\geq4$. If $u_1,u_2\in C((0,T];\widetilde{L}_{\sigma}^{d,\infty}(\mathbb{R}^d))\cap L^\beta(0,T;\widetilde{L}_{\sigma}^{d,\infty}(\mathbb{R}^d))$ are two mild solutions of $(1)$ in $(0,T)\times\mathbb{R}^d$ with the same force $f$ and initial data $u_{0}$, then $u_1=u_2$ in $(0,T)\times\mathbb{R}^d$.\\[5pt]
(\romannumeral2) Let $d\geq3$, $\beta>2d/(d-2)$, $u_{0}\in X^{d,\infty}_{\sigma}(\mathbb{R}^d)$ and $f$ satisfy $t^{\alpha}f\in L^\infty(0,T;L^{d,\infty}(\mathbb{R}^d))$ for some $\alpha<1$. If $u_{1},u_{2}\in C((0,T];\widetilde{L}_{\sigma}^{d,\infty}(\mathbb{R}^d))\cap L^\beta(0,T;\widetilde{L}_{\sigma}^{d,\infty}(\mathbb{R}^d))$ are two mild solutions of $(1)$ in $(0,T)\times\mathbb{R}^d$ with the same force $f$ and initial data $u_{0}$, then $u_{1}=u_{2}$ in $(0,T)\times\mathbb{R}^d$.\\[5pt]
\textbf{Remark 4.} It should be emphasized that Corollary 1 is the first result of the unconditional uniqueness without assuming a priori compactness nor smallness of the solution orbit $\{u(t)\}_{t\in[0,T]}$ in the scale critical space. We also obtain the following corollary directly from Theorem 1.\\[5pt]
\textbf{Corollary 2.} Let $d\geq3$, $\beta>2d/(d-2)$, $u_{0}\in \widetilde{L}^{d,\infty}_{\sigma}(\mathbb{R}^d)$, and $f\in L^1(0,T;L^{1}(\mathbb{R}^d))$ or $f\in L^1(0,T;L^{r,\infty}(\mathbb{R}^d))$ for some $r\in (1,\infty)$. Assume there are two distinct mild solutions $u_{1},u_{2}\in C((0,T];\widetilde{L}_{\sigma}^{d,\infty}(\mathbb{R}^d))\cap L^\beta(0,T;\widetilde{L}_{\sigma}^{d,\infty}(\mathbb{R}^d))$ of $(1)$ in $(0,T)\times\mathbb{R}^d$ with the same force $f$ and initial data $u_{0}$. Then for any $T_{0}>0$ there is no mild solution in $C([0,T_{0}];\widetilde{L}_{\sigma}^{d,\infty}({\mathbb{R}^d}))$ of (1) in $(0,T_{0})\times\mathbb{R}^d$ with the same $f$ and initial data $u_{0}$.\\[5pt]
\textbf{Remark 5.} (\romannumeral1) For any $T_{0}>0$, there is no mild solution in $C([0,T_{0}];\widetilde{L}_{\sigma}^{3,\infty}({\mathbb{R}^3}))$ of (1) in $(0,T_{0})\times\mathbb{R}^3$ with $f$ constructed in \cite{albritton2022non} and $u_{0}=0$.\\[5pt]
(\romannumeral2) The force $f$ constructed in \cite{albritton2022non} satisfies $f(t,x)=t^{-3/2}F(x/\sqrt{t})$, where $F(x)$ is smooth and compact in $\mathbb{R}^3$. Compared to Corollary 1, we notice that $\|f(t)\|_{L^{3,\infty}(\mathbb{R}^3)}\sim 1/t$, and $f\in BC((0,T];L^1(\mathbb{R}^3))$ but is not continuous at $t=0$ in $L^1(\mathbb{R}^3)$, which may explain why the uniqueness fails in $BC((0,T];L^3(\mathbb{R}^3))$.\vspace{5pt}

In the proof of Theorem 1, we may assume $T_{0}=T$ without loss of generality. Indeed, if we can prove $u_{1}=v$ in $(0,T_{0})\times\mathbb{R}^d$, then we have $u_{1}\in C([0,T];\widetilde{L}_{\sigma}^{d,\infty}(\mathbb{R}^d))$. Therefore, we can iterate the result to establish the uniqueness.

The strategy to prove Theorem 1 is sketched as follows. As in \cite{lions1998unicite}, we consider the dual problem of the Navier-Stokes equations:
\begin{equation}\begin{cases}
    -\partial_{t}\phi-(u\cdot\nabla)\phi-\Delta\phi-(\nabla\otimes\phi)\cdot v+\nabla\pi=F,\\
    \nabla\cdot\phi=0,\quad \phi(T)=0,
\end{cases}\end{equation}
where $u\in C((0,T];\widetilde{L}_{\sigma}^{d,\infty}(\mathbb{R}^d))\cap L^\beta(0,T;\widetilde{L}_{\sigma}^{d,\infty}(\mathbb{R}^d))$ with $\beta>2d/(d-2)$ and 
$v\in C([0,T];\widetilde{L}_{\sigma}^{d,\infty}(\mathbb{R}^d))$. The key idea is to derive the identity \begin{equation}\int_{\eta}^{T}\int_{\mathbb{R}^d}(u-v)\cdot F dxdt=\int_{\mathbb{R}^d}(u-v)(\eta,x)\cdot\phi(\eta,x)dx
\end{equation} for any $\eta\in(0,T)$, where $\phi$ is a solution to $(4)$. We note that, in \cite{lions1998unicite}, the identity $(5)$ is derived directly for $\eta=0$ under the assumption $u,v\in C([0,T];L^d(\mathbb{R}^d))$. In this paper, due to the lack of the condition $u\in C([0,T];\widetilde{L}_{\sigma}^{d,\infty}(\mathbb{R}^d))$, we first take positive $\eta$ and investigate the limit $\eta\to0$ in $(5)$. A crucial point is that, we can derive the regularity such as $\phi\in(L^{\infty}(0,T;L^2(\mathbb{R}^d)))^d$ by using the condition $v\in C([0,T];\widetilde{L}_{\sigma}^{d,\infty}(\mathbb{R}^d))$. Then we obtain $\int_{0}^{T}\int_{\mathbb{R}^d}(u-v)\cdot F dxdt=0$ from $(5)$ by showing $\lim_{\eta\to 0}\|(u-v)(\eta)\|_{L^2_x}=0$, which implies that $u-v=0$ in $(0,T)\times\mathbb{R}^d$. The key observation is that the $L^2$ norm is subcritical in view of scaling when $d\geq3$, and therefore, it is not difficult to derive $\lim_{\eta\to 0}\|(u-v)(\eta)\|_{L^2_x}=0$.

This paper is organized as follows. In Section 2 we introduce some well-known propositions. In Section 3 we give the proof of Theorem 1 and Corollary 1.
\section{Preliminaries}
In this section we recall some well-known results:\\
\textbf{Proposition 1} (\cite[Chapter 11, Theorem 11.2]{lemarie2002recent})\textbf{.} Let $f\in L^1(0,T;L^1(\mathbb{R}^d))$ or $f\in L^1(0,T;L^{r,\infty}(\mathbb{R}^d))$ for some $r\in(1,\infty)$. For $u\in L^\beta(0,T;\widetilde{L}_{\sigma}^{d,\infty}(\mathbb{R}^d))$ with $\beta\geq2$, $u$ is a mild solution of $(1)$ if and only if $u$ is a very weak solution of $(1)$.\\[5pt]
\textbf{Proposition 2} \cite[Chapter 11]{lemarie2002recent}\textbf{.} Let $1\leq i,j\leq d$ and $t>0$. The operator $O_{i,j,t}=\frac{1}{\Delta}\partial_i\partial_je^{t\Delta}$ is a convolution operator $O_{i,j,t}f=t^{-d/2}K_{i,j}(\cdot/\sqrt{t})\ast f$. Here $K_{i,j}(x)$ satisfies 
$$|\partial_{x}^{\alpha}K_{i,j}(x)|\leq C_{\alpha}\frac{1}{(1+|x|)^{d+|\alpha|}},\quad\alpha\in \mathbb{N}^d,\,\, x\in \mathbb{R}^d.$$
\textbf{Proposition 3.} (\romannumeral1) Let $1<s<r<\infty$. For any $f\in L^{s}(\mathbb{R}^d)$, we have
$$\left\lVert e^{t\Delta}\mathbb{P}f\right\rVert_{L^{r,1}}\leq C_{d,r,s}t^{\frac{d}{2}(\frac{1}{r}-\frac{1}{s})}\left\lVert f\right\rVert_{L^{s,\infty}},\quad t>0,$$
with $C_{d,r,s}$ depending only on $d$, $r$, $s$.\\
(\romannumeral2) Let $1<r,q<\infty$ or $r=q=\infty$. For any $f\in L^{1}(\mathbb{R}^d)$,
$$\left\lVert e^{t\Delta}\mathbb{P}f\right\rVert_{L^{r,q}}\leq C'_{d,r,q}t^{\frac{d}{2}(\frac{1}{r}-1)}\left\lVert f\right\rVert_{L^{1}},\quad t>0,$$
with $C'_{d,r,q}$ depending only on $d$, $r$, $q$.\\
(\romannumeral3) Let $1< s<r<\infty $. For any $F\in L^{s,\infty}(\mathbb{R}^d)$, we have
$$\left\lVert e^{t\Delta}\mathbb{P}\nabla \cdot F\right\rVert_{L^{r,1}}\leq C''_{d,r,s}t^{-\frac{1}{2}+\frac{d}{2}(\frac{1}{r}-\frac{1}{s})}\left\lVert F\right\rVert_{L^{s,\infty}},\quad t>0,$$ 
with $C''_{d,r,s}$ depending only on $d$, $r$, $s$.\\
(\romannumeral4) Let $1<r<\infty $, $1<q\leq\infty $ or $r=q=1$ or $r=q=\infty $. For any $F\in L^{r,q}(\mathbb{R}^d)$, we have
$$\left\lVert e^{t\Delta}\mathbb{P}\nabla \cdot F\right\rVert_{L^{r,q}}\leq C'''_{d,r,q}t^{-\frac{1}{2}}\left\lVert F\right\rVert_{L^{r,q}},\quad t>0,$$ 
with $C'''_{d,r,q}$ depending only on $d$, $r$, $q$.\\
\textbf{Proof.} The estimate directrly follows from the pointwise estimate stated in Proposition 2 and the Young inequality for the convolution in Lorentz space.

\section{Proof of the uniqueness}
In this section we give the proof of Theorem 1. Let us start from the following lemma.\\
\textbf{Lemma 1.} Let $d\geq3$, $\beta>2$, and $u,v\in L^{\beta}(0,T;L^{d,\infty}({\mathbb{R}^d}))$. For any $ F\in C_{c}^{\infty}((0,T)\times\mathbb{R}^d)$, there exists a solution $\phi^{\epsilon}\in W^{1,2}(0,T;W^{k,2}(\mathbb{R}^d))$ for any $k\in\mathbb{N}$ to the following mollified dual problem:
\begin{equation}
\begin{cases}
    -\partial_{t}\phi^{\epsilon}-(u_{\epsilon}\cdot\nabla)\phi^{\epsilon}-\Delta\phi^{\epsilon}-(\nabla\otimes\phi^{\epsilon})\cdot v_{\epsilon}+\nabla p=F,\\
    \nabla\cdot\phi^{\epsilon}=0,\quad \phi^{\epsilon}(T)=0.
\end{cases}
\end{equation}
Here $$u_{\epsilon}(t,x)=\int_{\mathbb{R}^d}\gamma_{\epsilon}(x-y)u(t,y)dy,$$ $$v_{\epsilon}(t,x)=\int_{\mathbb{R}^d}\gamma_{\epsilon}(x-y)v(t,y)dy,$$ $$\gamma_{\epsilon}=\epsilon^{-d}\eta(x/\epsilon), \quad \gamma\in C_{c}^{\infty}(B(0,1)),\quad\gamma\geq0\quad and\quad  \lVert\gamma\rVert_{L^1_{x}}=1.$$ 
\textbf{Proof.} Set $\tilde{u}_{\epsilon}(t,x)=u_{\epsilon}(T-t,x)$, $\tilde{v}_{\epsilon}(t,x)=v_{\epsilon}(T-t,x)$, $\tilde{F}(t,x)=F(T-t,x)$, and define $$A(\tilde{F})(t)=\int_{0}^{t}e^{(t-s)\Delta}\mathbb{P}\tilde{F}ds,$$
$$B(g)(t)=\int_{0}^{t}e^{(t-s)\Delta}\mathbb{P}((\tilde{u}_{\epsilon}\cdot\nabla)g+(\nabla\otimes g)\cdot \tilde{v}_{\epsilon})ds.$$
We can easily check from Proposition 3 that
$$\lVert A(\tilde{F})\rVert_{L^{\infty}(0,T;W^{k,p}(\mathbb{R}^d))}\leq C_{p}T\lVert \tilde{F}\rVert_{L^{\infty}(0,T;W^{k,p}(\mathbb{R}^d))}, \,\,k\in\mathbb{N},\,\,p\in(1,\infty),$$ 
and 
$$\lVert A(\tilde{F})\rVert_{L^{\infty}(0,T;W^{k,\infty}(\mathbb{R}^d))}\leq C\sqrt{T}\lVert \tilde{F}\rVert_{L^{\infty}(0,T;W^{k,d}(\mathbb{R}^d))},\,\,k\in\mathbb{N},$$
with $C_{p}$, $C$ independent of $\tilde{F}$, $k$ and $T$. 

Next we set $X_{t}=L^{\infty}(0,T;W^{1,2}(\mathbb{R}^d))\cap L^{\infty}(0,T;W^{1,\infty}(\mathbb{R}^d)).$ Let $t<T$ and $g\in (X_{t})^d$. By Proposition 3, we get 
$$\begin{aligned}
\lVert B(g)\rVert_{L^{\infty}_{t}(0,t;L_{x}^\infty)}&\leq C\Big\|\int_{0}^{t}(t-s)^{-1/2}\lVert(\tilde{u}_{\epsilon}\cdot\nabla)g+(\nabla\otimes g)\cdot \tilde{v}_{\epsilon}\rVert_{L_{x}^{d,\infty}}ds\Big\|_{L^{\infty}(0,t)}\\
&\leq C \|t^{-1/2}\|_{L^{\beta/(\beta-1)}(0,t)}(\lVert \tilde{u}\rVert_{L_{t}^{\beta}(0,t;L_{x}^{d,\infty})}+\lVert \tilde{v}\rVert_{L_{t}^{\beta}(0,t;L_{x}^{d,\infty})})\lVert \nabla g\rVert_{L^{\infty}_{t}(0,t;L_{x}^\infty)}\\
&\leq C t^{1/2-1/\beta}(\lVert \tilde{u}\rVert_{L_{t}^{\beta}(0,T;L_{x}^{d,\infty})}+\lVert \tilde{v}\rVert_{L_{t}^{\beta}(0,T;L_{x}^{d,\infty})})\lVert \nabla g\rVert_{L^{\infty}_{t}(0,t;L_{x}^\infty)},
\end{aligned}$$
$$\begin{aligned}
\lVert B(g)\rVert_{L^{\infty}_{t}(0,t;L_{x}^2)}&\leq C\Big\|\int_{0}^{t}(t-s)^{-1/2}
\lVert (\tilde{u}_{\epsilon}\cdot\nabla)g+(\nabla\otimes g)\cdot \tilde{v}_{\epsilon}\rVert_{L_{x}^{2d/(2+d),\infty}}ds\Big\|_{L^{\infty}(0,t)}\\
&\leq C\|t^{-1/2}\|_{L^{\beta/(\beta-1)}(0,t)}(\lVert \tilde{u}\rVert_{L_{t}^{\beta}(0,t;L_{x}^{d,\infty})}+\lVert \tilde{v}\rVert_{L_{t}^{\beta}(0,t;L_{x}^{d,\infty})})\lVert \nabla g\rVert_{L^{\infty}_{t}(0,t;L_{x}^2)}\\
&\leq C t^{1/2-1/\beta}(\lVert \tilde{u}\rVert_{L_{t}^{\beta}(0,T;L_{x}^{d,\infty})}+\lVert \tilde{v}\rVert_{L_{t}^{\beta}(0,T;L_{x}^{d,\infty})})\lVert \nabla g\rVert_{L^{\infty}_{t}(0,t;L_{x}^2)},
\end{aligned}$$
$$\begin{aligned}
\lVert \nabla B(g)\rVert_{L^{\infty}_{t}((0,t);L_{x}^\infty)}&\leq C\Big\|\int_{0}^{t}(t-s)^{-1/2}
\lVert (\tilde{u}_{\epsilon}\cdot\nabla)g+(\nabla\otimes g)\cdot \tilde{v}_{\epsilon}\rVert_{L_{x}^{\infty}}ds\Big\|_{L^{\infty}(0,t)}\\
&\leq C\|t^{-1/2}\|_{L^{\beta/(\beta-1)}(0,t)}(\|\tilde{u}_{\epsilon}\|_{L^{\beta}(0,t;L^{\infty}_{x})}+\|\tilde{v}_{\epsilon}\|_{L^{\beta}(0,t;L^{\infty}_{x})}) \lVert\nabla g\rVert_{L^{\infty}_{t}(0,t;L_{x}^\infty)}\\
&\leq C\epsilon^{-1}t^{1/2-1/\beta}(\|\tilde{u}\|_{L^{\beta}(0,T;L^{d,\infty}_{x})}+\|\tilde{v}\|_{L^{\beta}(0,T;L^{d,\infty}_{x})})\lVert \nabla g\rVert_{L^{\infty}_{t}(0,t;L_{x}^\infty)},
\end{aligned}$$
$$\begin{aligned}
\lVert \nabla B(g)\rVert_{L^{\infty}_{t}(0,t;L_{x}^2)}&\leq C\Big\|\int_{0}^{t}(t-s)^{-1/2}
\lVert (\tilde{u}_{\epsilon}\cdot\nabla)g+(\nabla\otimes g)\cdot \tilde{v}_{\epsilon}\rVert_{L_{x}^{2}}ds\Big\|_{L^{\infty}(0,t)}\\
&\leq C\|t^{-1/2}\|_{L^{\beta/(\beta-1)}(0,t)}(\|\tilde{u}_{\epsilon}\|_{L^{\beta}(0,t;L^{\infty}_{x})}+\|\tilde{v}_{\epsilon}\|_{L^{\beta}(0,t;L^{\infty}_{x})}) \lVert \nabla g\rVert_{L^{\infty}_{t}(0,t;L_{x}^2)}\\
&\leq C\epsilon^{-1}t^{1/2-1/\beta}(\|\tilde{u}\|_{L^{\beta}(0,T;L^{d,\infty}_{x})}+\|\tilde{v}\|_{L^{\beta}(0,T;L^{d,\infty}_{x})})\lVert \nabla g\rVert_{L^{\infty}_{t}(0,t;L_{x}^2)}.
\end{aligned}$$
Hence, we have $\lVert B(g)\rVert_{X_{t}}\leq C(\epsilon)t^{1/2-1/\beta}\lVert g\rVert_{X_{t}}$ where $C(\epsilon)$ is independent of $g$ and $t$. Set $g^{(0)}=A(\tilde{F})\in (X_{t})^d$, $g^{(n)}=g^{(0)}+B(g^{(n-1)})$ for $n\geq 1$, we get 
$$\lVert g^{(1)}-g^{(0)}\rVert_{X_{t}}\leq C(\epsilon)t^{1/2-1/\beta}\lVert g^{(0)}\rVert_{X_{t}},$$
$$\lVert g^{(n+1)}-g^{(n)}\rVert_{X_{t}}\leq C(\epsilon)t^{1/2-1/\beta}\lVert g^{(n)}-g^{(n-1)}\rVert_{X_{t}},\quad n\geq1.$$
Therefore, for any fixed $\epsilon>0$ and for $\delta$ small enough such that $C(\epsilon)t^{1/2-1/\beta}<1$, by fixed point theorem we can construct a solution $g^{\epsilon}\in X_{\delta}$ satisfying
$g^{\epsilon}=A(\tilde{F})+B(g^{\epsilon})$. Moreover, for any $t_{0}\in(0,\delta)$, by regarding $g(t_{0})$ as the initial data, we can extend $g^{\epsilon}$ from $t_{0}$ to $t_{0}+\delta$ so that $g^{\epsilon}=A(\tilde{F})+B(g^{\epsilon})$ holds. Indeed, let us set $g^{(0)}(t)=e^{t\Delta}g^{\epsilon}(t_{0})+\int_{0}^{t}e^{(t-s)\Delta}\mathbb{P}\tilde{F}(t_{0}+s)ds$ and $g^{(n)}=g^{(0)}+B(g^{(n-1)})$ for $n\geq 1$. Since the choice of $\delta$ satisfying
$C(\epsilon)t^{1/2-1/\beta}<1$ is independent of $t_{0}$, we obtain $g^{\epsilon}\in X_{T}$ such that
\begin{equation}
    g^{\epsilon}(t)=\int_{0}^{t}e^{(t-s)\Delta}\mathbb{P}\tilde{F}ds+\int_{0}^{t}e^{(t-s)\Delta}\mathbb{P}((\tilde{u}_{\epsilon}\cdot\nabla)g^{\epsilon}+(\nabla\otimes g^{\epsilon})\cdot \tilde{v}_{\epsilon})ds,\quad t\in[0,T].
\end{equation}

We notice that $\phi^{\epsilon}(t)=g^{\epsilon}(T-t)$ is the solution of $(6)$ in the sense of distributions. Now we prove $\phi^{\epsilon}\in W^{1,2}(0,T;W^{k,2}(\mathbb{R}^d))$. Assume $g^{\epsilon}\in (L^{\infty}(0,T;W^{k,2}(\mathbb{R}^d)))^d$ for some $k\in\mathbb{N}$, then 
$$\begin{aligned}
   &\lVert (-\Delta)^{k/2+m}g^{\epsilon}\rVert_{L_{t}^{\infty}L_{x}^2}\\
   &\leq \lVert(-\Delta)^{k/2+m}A(\tilde{F})\rVert_{L_{t}^{\infty}L_{x}^2}+\big\lVert\int_{0}^{t}(-\Delta)^{1/2+m}e^{(t-s)\Delta}\mathbb{P}(-\Delta)^{(k-1)/2} ((\tilde{u}_{\epsilon}\cdot\nabla)g^{\epsilon}+(\nabla\otimes g^{\epsilon})\cdot \tilde{v}_{\epsilon})ds\big\rVert_{L_{t}^{\infty}L_{x}^2}\\
   &\leq C_{1}T+C_{2}(\epsilon)\|t^{-1/2-m}\|_{L^{\beta/(\beta-1)}(0,T)}(\|\tilde{u}\|_{L^{\beta}(0,T;L^{d,\infty})}+\|\tilde{v}\|_{L^{\beta}(0,T;L^{d,\infty})})\lVert g^{\epsilon}\rVert_{L_{t}^{\infty}W_{x}^{k,2}}\\
   &\leq  C_{1}T+C_{2}(\epsilon)T^{1/2-1/\beta-m}(\|\tilde{u}\|_{L^{\beta}(0,T;L^{d,\infty})}+\|\tilde{v}\|_{L^{\beta}(0,T;L^{d,\infty})})\lVert g^{\epsilon}\rVert_{L_{t}^{\infty}W_{x}^{k,2}},
\end{aligned}$$
where $m<1/2-1/\beta$.
By induction we obtain $g^{\epsilon}\in L^{\infty}(0,T;W^{k,2}(\mathbb{R}^d))$ for any $k\in\mathbb{N}$. As for $\partial_{t}g^{\epsilon}$, going back to $(6)$, we obtain that
$$\lVert\partial_{t}g^{\epsilon}\rVert_{L^{2}(0,T;W^{k,2}(\mathbb{R}^d))}\leq C(\epsilon,T)(
\lVert g^{\epsilon}\rVert_{L^{\infty}(0,T;W^{k+2,2}(\mathbb{R}^d))}+\lVert \tilde{F}\rVert_{L^{\infty}(0,T;W^{k,2}(\mathbb{R}^d))}).$$
This completes the proof. $\hfill\square$\\[5pt]
\textbf{Proof of Theorem 1.} Let $u_{1}\in C((0,T];\widetilde{L}^{d,\infty}_{\sigma}(\mathbb{R}^d))\cap L^\beta(0,T;\widetilde{L}^{d,\infty}_{\sigma}(\mathbb{R}^d))$ with $\beta>2d/(d-2)$, $v\in C([0,T_{0}];\widetilde{L}_{\sigma}^{d,\infty}(\mathbb{R}^d))$ be the mild solutions assumed in Theorem 1. As stated in the introduction, we may assume $T_{0}=T$ without loss of generality. Set $u_{1}=u$. Our purpose is to prove $u=v$ on $(0,T)\times\mathbb{R}^d$. According to Proposition 1, the mild solution of $(1)$ is the very weak solution of $(1)$. Set $w=u-v$. For any $\eta\in(0,T)$, $\phi\in
C_{c,\sigma}^{\infty}([\eta,T]\times\mathbb{R}^d)$ with $\phi(T)=0$, we have \begin{equation}
\int_{\mathbb{R}^d}w(\eta,x)\cdot\phi(\eta,x)dx+\int_{\eta}^{T}\int_{\mathbb{R}^d}w\cdot(\partial_{t}\phi+(u\cdot\nabla)\phi+\Delta\phi+(\nabla\otimes\phi)\cdot v)dxdt=0.
\end{equation}
For any $F\in C_{c}^{\infty}((0,T)\times\mathbb{R}^d)$, we consider the dual problem of $(8)$:
\begin{equation}
\begin{cases}
    -\partial_{t}\phi-(u\cdot\nabla)\phi-\Delta\phi-(\nabla\otimes\phi)\cdot v+\nabla\pi=F,\\
    \nabla\cdot\phi=0,\quad \phi(T)=0.
\end{cases}
\end{equation}
If there is a solution $\phi$ to $(9)$ in a suitable Sobolev space, then formally we have $\int_{\mathbb{R}^d}w(\eta,x)\cdot\phi(\eta,x)dx-\int_{\eta}^{T}\int_{\mathbb{R}^d}w\cdot Fdxdt=0$ for $\eta\in(0,T]$. To ensure such $\phi$ exists, we start from the mollified equations stated in Lemma 1, for which we have constructed $\phi^{\epsilon}\in W^{1,2}(0,T;W^{k,2}(\mathbb{R}^d))$. Next we prove that for any fixed $\eta\in(0,T)$, $p\in(1,d)$, and $p^*\in(1,\infty)$, $\phi^{\epsilon}$ is uniformly bounded with respect to $\epsilon$ in 
$L^{\infty}(0,T;L^2(\mathbb{R}^d))\cap L^{p^*}(\eta,T;W^{2,p}(\mathbb{R}^d))\cap W^{1,p^*}(\eta,T;L^{p}(\mathbb{R}^d))$. For this purpose, we start from the energy estimate for the terminal value problem (9). For any $t\in[0,T]$, we have
$$\begin{aligned}
    \frac{1}{2}\partial_{t}\int_{\mathbb{R}^d}|\phi^{\epsilon}|^2dx-\int_{\mathbb{R}^d}|\nabla\phi^{\epsilon}|^2dx&=-\int_{\mathbb{R}^d}\phi^{\epsilon}\cdot(\nabla\otimes\phi^{\epsilon})\cdot v_{\epsilon}dx-\int_{\mathbb{R}^d}F\cdot\phi^{\epsilon} dx
\end{aligned}$$
Here we have used $\int_{\mathbb{R}^d}\phi^{\epsilon}\cdot(u_{\epsilon}\cdot\nabla)\phi^{\epsilon}dx=0$ because $u_{\epsilon}$ is divergence free.
Since $v\in C([0,T];\widetilde{L}_{\sigma}^{d,\infty})$, for any $\delta>0$ we can decompose $v$ into $v=v^{\delta}+\bar{v}^{\delta}$ such that $\lVert v^{\delta}\rVert_{C([0,T];L^{d,\infty})}\leq \delta$ and $\|\bar{v}^{\delta}\|_{L^{\infty}_{t,x}}\leq K(\delta)$ (see \cite{okabe2021remark}).
Therefore, $v_{\epsilon}=v^{\delta}_{\epsilon}+\bar{v}^{\delta}_{\epsilon}$ with $v^{\delta}_{\epsilon}=v^{\delta}\ast\eta_{\epsilon}$ and $\bar{v}^{\delta}_{\epsilon}=\bar{v}^{\delta}\ast\eta_{\epsilon}$. We easily check that $\lVert v^{\delta}_{\epsilon}\rVert_{C([0,T];L^d)}\leq \lVert v^{\delta}\rVert_{C([0,T];L^{d,\infty})}\leq\delta$ and $\|\bar{v}_{\epsilon}^{\delta}\|_{L^{\infty}_{t,x}}\leq\|\bar{v}^{\delta}\|_{L^{\infty}_{t,x}}\leq K(\delta)$. Using this decomposition, we derive
 \begin{align*}
    &\frac{1}{2}\partial_{t}\int_{\mathbb{R}^d}|\phi^{\epsilon}|^2dx-\int_{\mathbb{R}^d}|\nabla\phi^{\epsilon}|^2dx\geq -\int_{\mathbb{R}^d}(|v^{\delta}_{\epsilon}|+|\bar{v}^{\delta}_{\epsilon}|)|\nabla\phi^{\epsilon}||\phi^{\epsilon}|dx-\int_{\mathbb{R}^d}|F||\phi^{\epsilon}|dx\\
    &\geq-\int_{\mathbb{R}^d}|v^{\delta}_{\epsilon}||\nabla\phi^{\epsilon}||\phi^{\epsilon}|dx-\int_{\mathbb{R}^d}K(\delta)|\nabla\phi^{\epsilon}||\phi^{\epsilon}|dx-\frac{1}{2}\int_{\mathbb{R}^d}|F|^2dx-\frac{1}{2}\int_{\mathbb{R}^d}|\phi^{\epsilon}|^2dx.
\end{align*}
By using the Hölder inequality and the Young inequality in the last line, we obtain 
\begin{align*}
&\frac{1}{2}\partial_{t}\int_{\mathbb{R}^d}|\phi^{\epsilon}|^2dx-\int_{\mathbb{R}^d}|\nabla\phi^{\epsilon}|^2dx\\
    &\geq -\lVert v^{\delta}_{\epsilon}\rVert_{L^{d,\infty}}\lVert \nabla\phi^{\epsilon}\rVert_{L^2}\lVert \phi^{\epsilon}\rVert_{L^{\frac{2d}{d-2},2}}-\frac{K^2(\delta)}{2}\int_{\mathbb{R}^d}|\phi^{\epsilon}|^2dx-\frac{1}{2}\int_{\mathbb{R}^d}|\nabla\phi^{\epsilon}|^2dx\\
    &-\frac{1}{2}\int_{\mathbb{R}^d}|F|^2dx-\frac{1}{2}\int_{\mathbb{R}^d}|\phi^{\epsilon}|^2dx\\
    &\geq -C\delta\lVert\nabla\phi^{\epsilon}\rVert^2_{L^2}-C(\delta)\int_{\mathbb{R}^d}|\phi^{\epsilon}|^2dx-\frac{1}{2}\int_{\mathbb{R}^d}|\nabla\phi^{\epsilon}|^2dx-\frac{1}{2}\int_{\mathbb{R}^d}|F|^2dx,
    \end{align*}
where $C$, $C(\delta)$ are independent of $\phi^{\epsilon}$ and $t$. We have used the Sobolev inequality $\lVert\phi^{\epsilon}\rVert_{L^{2d/(d-2),2}}\leq C\lVert\nabla\phi^{\epsilon}\rVert_{L^2}$ for $d\geq3$ in the last line. Take $\delta$ small enough such that $C\delta=1/4$. Then we obtain
\begin{equation}
    \frac{1}{2}\partial_{t}\int_{\mathbb{R}^d}|\phi^{\epsilon}|^2dx\geq \frac{1}{4}\int_{\mathbb{R}^d}|\nabla\phi^{\epsilon}|^2dx-C'\int_{\mathbb{R}^d}|\phi^{\epsilon}|^2dx-\frac{1}{2}\int_{\mathbb{R}^d}|F|^2dx,\quad t\in[0,T],
\end{equation}
where $C'$ is independent of $\phi^{\epsilon}$ and $t$. According to the Grönwall's inequality from the terminal time $T$, we get 
$$\sup_{t\in[0,T]}\lVert\phi^{\epsilon}(t)\rVert^2_{L_x^2}\leq \|F\|^2_{L^{\infty}_{t}L_{x}^2}Te^{2C'T}=C'',$$
where $C''$ is independent of $\phi^{\epsilon}$ and $t$. Hence $\phi^{\epsilon}$ is uniformly bounded in 
$L^{\infty}(0,T;L^2)$ with respect to $\epsilon$. Moreover, by integrating $(10)$ on $t$, we observe that $\phi^{\epsilon}$ is uniformly bounded in $L^2((0,T);H^1)$ with respect to $\epsilon$. Fix $\eta\in(0,T)$, we notice that $u\in L^{\infty}(\eta,T;L^{d,\infty})$. This yields that 
$|u_{\epsilon}||\nabla\phi^{\epsilon}|$ and $|v_{\epsilon}||\nabla\phi^{\epsilon}|$ are uniformly bounded in $L^{2}(\eta,T;L^{2d/(2+d),2})$ with respect to $\epsilon$. By going back to $(7)$, we get that 
$|\nabla\phi^{\epsilon}|$ is uniformly bounded in $L^q((\eta,T);L^{2d/(2+d),2})$ for any $q\in(2,\infty)$, since
$$
\begin{aligned}
 &\big\|\int_{0}^{t}\nabla e^{(t-s)\Delta}\mathbb{P}((\tilde{u}_{\epsilon}\cdot\nabla)g^{\epsilon}+(\nabla\otimes g^{\epsilon})\cdot\tilde{v}_{\epsilon})ds\big\|_{L^q(0,T-\eta;L^{\frac{2d}{2+d},2})}\\
 &\leq C\big\|\int_{0}^{t}(t-s)^{-\frac{1}{2}}\lVert(|\tilde{u}_{\epsilon}|+|\tilde{v}_{\epsilon}|)\nabla g^{\epsilon}(s)\rVert_{L^{\frac{2d}{2+d},2}}ds\big\|_{L^q(0,T-\eta)}\\
 &\leq C\lVert t^{-\frac{1}{2}}\rVert_{L^{\frac{2q}{q+2}}(0,T-\eta)}\lVert(|\tilde{u}_{\epsilon}|+|\tilde{v}_{\epsilon}|)\nabla g^{\epsilon}\rVert_{L^{2}(0,T-\eta;L^{\frac{2d}{2+d},2})}. 
\end{aligned}
$$
Again, we obtain that $|u_{\epsilon}||\nabla\phi^{\epsilon}|$ and $|v_{\epsilon}||\nabla\phi^{\epsilon}|$ are uniformly bounded in $L^{q}(\eta,T;L^{2d/(4+d),2})$ with $d\geq5$ or $L^{2}(\eta,T;L^{r,2})$ with $d=3$, $4$ and any $r\in (1,2d/(2+d))$. Iterating this procedure, we deduce that $|\nabla\phi^{\epsilon}|$ is uniformly bounded with respect to $\epsilon$ in $L^q((\eta,T);L^{r})$ for any $q\in[1,\infty)$, $r\in(1,2d/(2+d))$.

Now we can prove that, $\phi^{\epsilon}$ is uniformly bounded in 
$L^{p^*}(\eta,T;W^{2,p})\cap W^{1,p^{*}}(\eta,T;L^{p})$ for any $p\in(1,d)$, $p^*\in(1,\infty)$, and $\eta\in(0,T)$. According the $L_{t}^{p^*}L_{x}^{p}$ regularity of the heat kernel, we have 
\begin{align}
&\big\|\partial_{t}\phi^{\epsilon}\big\|_{L^{p^*}(\eta,T;L^{p})}+\big\|\Delta\phi^{\epsilon}\big\|_{L^{p^*}(\eta,T;L^{p})}
+\big\|\nabla\pi\big\|_{L^{p^*}(\eta,T;L^{p})}\\\notag
&\leq C\big\|F\big\|_{L^{p^*}(\eta,T;L^{p})}+C\big\||u_{\epsilon}||\nabla\phi^{\epsilon}|\big\|_{L^{p^*}(\eta,T;L^{p})}+C\big\||v_{\epsilon}||\nabla\phi^{\epsilon}|\big\|_{L^{p^*}(\eta,T;L^{p})}.
\end{align}

Since $u,v\in C([\eta,T];\widetilde{L}_{\sigma}^{d,\infty}(\mathbb{R}^d))$, we decompose $u=u^{\delta}+\bar{u}^{\delta}$, $v=v^{\delta}+\bar{v}^{\delta}$ on $[\eta,T]$, so that $\|u^{\delta}\|_{C([\eta,T];\widetilde{L}_{\sigma}^{d,\infty}(\mathbb{R}^d))} \leq \delta $, $\|\bar{u}^{\delta}\|_{L^{\infty}((\eta,T)\times\mathbb{R}^d)}\leq K(\delta)$, $\|v^{\delta}\|_{C([\eta,T];\widetilde{L}_{\sigma}^{d,\infty}(\mathbb{R}^d))}\leq \delta $, and $\|\bar{v}^{\delta}\|_{L^{\infty}((\eta,T)\times\mathbb{R}^d)}\leq K(\delta)$. We notice that $u_{\epsilon}=u^{\delta}_{\epsilon}+\bar{u}^{\delta}_{\epsilon}$ and $v_{\epsilon}=v^{\delta}_{\epsilon}+\bar{v}^{\delta}_{\epsilon}$ on $[\eta,T]$, satisfying $\|u^{\delta}_{\epsilon}\|_{C([\eta,T];\widetilde{L}_{\sigma}^{d,\infty}(\mathbb{R}^d))} \leq \delta $, $\|\bar{u}^{\delta}_{\epsilon}\|_{L^{\infty}((\eta,T)\times\mathbb{R}^d)}\leq K(\delta)$, $\|v^{\delta}_{\epsilon}\|_{C([\eta,T];\widetilde{L}_{\sigma}^{d,\infty}(\mathbb{R}^d))} \leq \delta $, and $\|\bar{v}^{\delta}_{\epsilon}\|_{L^{\infty}((\eta,T)\times\mathbb{R}^d)}\leq K(\delta)$. Then,
\begin{align}
    \big\||u_{\epsilon}||\nabla\phi^{\epsilon}|\big\|_{L^{p^*}(\eta,T;L^{p})}+\big\||v_{\epsilon}||\nabla\phi^{\epsilon}|\big\|_{L^{p^*}(\eta,T;L^{p})}&\leq 2K(\delta)\big\|\nabla\phi^{\epsilon}\big\|_{L^{p^*}(\eta,T;L^{p})}\notag\\ 
    &+2\delta\big\|\nabla\phi^{\epsilon}\big\|_{L^{p^*}(\eta,T;L^{\frac{dp}{d-p},p})}.
\end{align}
Here the constant $K(\delta)$ may depend on $\delta$ and $\eta$. Applying the Sobolev and the Hölder inequalities, we deduce that 
\begin{align*}
\|\nabla\phi^{\epsilon}\|_{L^{p^*}(\eta,T;L^{\frac{dp}{d-p},p})}&\leq 
C\|\Delta\phi^{\epsilon}\|_{L^{p^*}(\eta,T;L^p)}\\
&\leq C\|F\|_{L^{p^*}(\eta,T;L^{p})}+CK(\delta)\|\nabla\phi^{\epsilon}\|_{L^{p^*}(\eta,T;L^{p})}+C\delta\|\nabla\phi^{\epsilon}\|_{L^{p^*}(\eta,T;L^{\frac{dp}{d-p},p})}.
\end{align*}
Take $\delta>0$ small enough, we have\begin{align*}
    \|\nabla\phi^{\epsilon}\|_{L^{p^*}(\eta,T;L^{\frac{dp}{d-p},p})}&\leq
C\|F\|_{L^{p^*}(\eta,T;L^{p})}+C\|\nabla\phi^{\epsilon}\|_{L^{p^*}(\eta,T;L^{p})}\\
&\leq C\|F\|_{L^{p^*}(\eta,T;L^{p})}+C
\|\nabla\phi^{\epsilon}\|^{\theta}_{L^{p^*}(\eta,T;L^{\frac{dp}{d-p},p})}\|\nabla\phi^{\epsilon}\|_{L^{p^*}((\eta,T);L^{r})}^{1-\theta},
\end{align*}
where $\theta(1/p-1/d)+(1-\theta)/r=1/p$. Since $\nabla\phi^\epsilon$ is uniformly bounded in $L^q(\eta,T;L^r)$ for $q\in [1,\infty)$ and $r\in (1,2d/(2+d))$, we obtain that $\nabla\phi^{\epsilon}$ is uniformly bounded in $L^{p^*}(\eta,T;L^{\frac{dp}{d-p},p})$. By interpolation, we obtain that $\nabla\phi^{\epsilon}$ is uniformly bounded in 
$L^{p^*}(\eta,T;L^{p})$ with respect to $\epsilon$, since $1<p<dp/(d-p)$. By going back to $(11)$ and $(12)$, we obtain that $\phi^{\epsilon}$ is uniformly bounded in 
$L^{p^*}(\eta,T;W^{2,p})\cap W^{1,{p^*}}(\eta,T;L^{p})$. Hence, we can take a subsequence $\{\phi^{\epsilon}\}_{\epsilon}$ so that for any $\eta\in(0,T)$, ${p^*}\in(1,\infty)$, $p\in(1,d)$, $\{\phi^{\epsilon}\}_{\epsilon}$ $*$-weakly converges to $\phi$ in $L^{\infty}(0,T;L^2(\mathbb{R}^d))\cap L^{p^*}(\eta,T;W^{2,p}(\mathbb{R}^d))\cap W^{1,{p^*}}(\eta,T;L^{p}(\mathbb{R}^d))$.  Let us recall that $\phi^{\epsilon}$ is the solution of $(6)$ in Lemma 1. Then for any fixed $\eta\in(0,T)$ and $h\in C_{c,\sigma}^{\infty}([\eta,T]\times\mathbb{R}^d)$ with $h(\eta,x)=0$, we have
\begin{align*}
&\Big|\int_{\eta}^{T}\int_{\mathbb{R}^d}h\cdot(\partial_{t}\phi+(u\cdot\nabla)\phi+\Delta\phi+(\nabla\otimes\phi)\cdot v)dxdt+\int_{\eta}^{T}\int_{\mathbb{R}^d}h(t,x)\cdot F(t,x)dxdt\Big|\\
\leq&\varlimsup_{\epsilon\to 0}\Big|\int_{\eta}^{T}\int_{\mathbb{R}^d}h\cdot(\partial_{t}(\phi-\phi^{\epsilon})+
(u\cdot\nabla)(\phi-\phi^{\epsilon})+\Delta(\phi-\phi^{\epsilon})+(\nabla\otimes(\phi-\phi^{\epsilon}))\cdot v)dxdt\Big|\\
+&\varlimsup_{\epsilon\to 0}\Big|\int_{\eta}^{T}\int_{\mathbb{R}^d}h\cdot((
(u-u_{\epsilon})\cdot\nabla)\phi^{\epsilon}+(\nabla\otimes\phi^{\epsilon})\cdot(v-v_{\epsilon}))dxdt\Big|\\
\leq&\varlimsup_{\epsilon\to 0}\Big|\int_{\eta}^{T}\int_{\mathbb{R}^d}h\cdot((
(u-u_{\epsilon})\cdot\nabla)\phi^{\epsilon}+(\nabla\otimes\phi^{\epsilon})\cdot(v-v_{\epsilon}))dxdt\Big|.
\end{align*}
We recall that $\widetilde{L}^{d,\infty}=\overline{L^{d,\infty}\cap L^{s}}^{L^{d,\infty}}$ for any $s\in (d,\infty)$, so $u$ can be decomposed into $u=u^{\delta}+\bar{u}^{\delta}$ s.t. $u^{\delta}\in C([\eta,T];L^{d,\infty})$ with norm smaller than $\delta$, and $\bar{u}^{\delta}\in C([\eta,T];L^{s})$ (see \cite{okabe2021remark}). Using the uniform continuity of $\bar{u}^{\delta}$ and $\bar{u}^{\delta}_{\epsilon}$ in $C([\eta,T];L^s)$, one can easily check that $\bar{u}^{\delta}_{\epsilon}$ converges to $\bar{u}^{\delta}$ in $C([\eta,T];L^{s})$ as $\delta$ fixed and $\epsilon$ goes to 0. Therefore, $u-u_{\epsilon}$ can be docomposed into $u^{\delta}-u^{\delta}_{\epsilon}$ and $\bar{u}^{\delta}-\bar{u}^{\delta}_{\epsilon}$ with $\|u^{\delta}-u^{\delta}_{\epsilon}\|_{C([\eta,T];L^{d,\infty})}\leq 2\delta$ and $\|\bar{u}^{\delta}-\bar{u}^{\delta}_{\epsilon}\|_{C([\eta,T];L^{s})}\leq \delta$ for $\epsilon$ very small. This implies $$\lim_{\epsilon\to 0}\int_{\eta}^{T}\int_{\mathbb{R}^d}h\cdot((
(u-u_{\epsilon})\cdot\nabla)\phi^{\epsilon}+(\nabla\otimes\phi^{\epsilon})\cdot(v-v_{\epsilon}))dxdt=0.$$

Therefore, $\phi$ is the solution of $(9)$. Although $\phi$ may not be compactly supported nor smooth, one can verify that $\phi$ can be taken as the test function in $(8)$. Indeed, by the density argument, there exists $\phi_{n}\in C_{c,\sigma}^{\infty}([\eta,T]\times\mathbb{R}^d)$, $\phi_{n}(T)=0$, which converges to $\phi$ in $L^{p^*}(\eta,T;W^{2,p}(\mathbb{R}^d))\cap W^{1,p^{*}}(\eta,T;L^{p}(\mathbb{R}^d))$. Particularly, $\phi_{n}(\eta,x)$ converges to $\phi(\eta,x)$ in $L^{p}(\mathbb{R}^d)$. Since $\phi$ is the solution of (9), we go back to $(8)$ and obtain that
\begin{equation}
\int_{\mathbb{R}^d}w(\eta,x)\cdot\phi(\eta,x)dx=\int_{\eta}^{T}\int_{\mathbb{R}^d}w\cdot Fdxdt,\quad\eta\in(0,T).
\end{equation}
Besides, we have the following estimate of $w=u-v$ : 
$$\begin{aligned}
\left\lVert w(t,\cdot)\right\rVert_{L^{s_1}(0,\eta;L^{d/2,\infty})}&= 
\Big\|\int_{0}^{t}e^{(t-s)\Delta}\mathbb{P}\nabla\cdot((w\otimes u)+(v\otimes w))ds\Big\|_{L^{s_1}(0,\eta;L^{d/2,\infty})}\\
&\leq C\Big\|\int_{0}^{t}(t-s)^{-1/2}\left\lVert w(s,\cdot)\right\rVert_{L^{d,\infty}}(\left\lVert u(s,\cdot)\right\rVert_{L^{d,\infty}}+\left\lVert v(s,\cdot)\right\rVert_{L^{d,\infty}})ds\Big\|_{L^{s_1}(0,\eta)}\\
&\leq C\|t^{-1/2}\|_{L^{r_1}(0,\eta)}\left\lVert w\right\rVert_{L^\beta(0,\eta;L^{d,\infty})}(\left\lVert u\right\rVert_{L^\beta(0,\eta;L^{d,\infty})}+\left\lVert v\right\rVert_{L^\beta(0,\eta;L^{d,\infty})})\\
&\leq C \eta^{1/{r_1}-1/2}\left\lVert w\right\rVert_{L^\beta(0,\eta;L^{d,\infty})}(\left\lVert u\right\rVert_{L^\beta(0,T;L^{d,\infty})}+\left\lVert v\right\rVert_{L^\beta(0,T;L^{d,\infty})}),
\end{aligned}$$
where $1+1/s_1=1/r_1+2/\beta$, $1/r_1>1/2$. We set $1/s_0=1/\beta$, by the real interpolation method, we can iterate the estimate above from $L^{s_{m-1}}L^{d/m}$ to $L^{s_m}L^{d/(m+1)}$ with $m\in(1,d-1)$, $1/s_m=1/s_{m-1}+1/\beta+1/r_m-1$, $1/r_m>1/2$. Therefore, there exists some $\gamma_m>0$ such that
$$\|w(t,\cdot)\|_{L^{s_m}(0,\eta;L^{d/(m+1)})}\leq C_m\eta^{\gamma_m},\quad m\in(1,d-1).$$ Taking $1/r_i=1/r=1-d/(d\beta-2\beta)>1/2$ for any $i\in[1,d/2-1]$, we obtain that $1/s_{d/2-1}=1/\beta+(d/2-1)(1/\beta+1/r-1)=0$. Hence
$$\|w(t,\cdot)\|_{L^{\infty}(0,\eta;L^2)}\leq C\eta^{1/r-1/2}.$$
This implies $\lim_{\eta\to0}\left\lVert w(\eta,\cdot)\right\rVert_{L^2}=0$ for $\beta>2d/(d-2)$. Consequently, $\int_{0}^{T}\int_{\mathbb{R}^d}w\cdot Fdxdt=0$ for any $F\in C_{c}^{\infty}((0,T)\times\mathbb{R}^d)$. This completes the proof.$\hfill\square$\\[5pt]
\textbf{Proof of Corollary 1.} (i) in Corollary 1 is a directly corollary of Theorem 1 and \cite[Theorem 2.3]{okabe2021remark}. As for (ii), the construction of a mild solution $u\in C([0,T_{0}];\widetilde{L}_{\sigma}^{d,\infty}(\mathbb{R}^d))$ for some $T_{0}>0$ with the initial data and force in (ii) is standard, which can be checked in \cite[Chpater 15, Theorem 15.3]{lemarie2002recent}, noticing that we only need to replace $L^d$ by $L^{d,\infty}$ and $L^{2d}$ by $L^{2d,\infty}$.

Define $\mathcal{E}_{T}=\{f\in L^{\infty}(0,T;L^{d,\infty}(\mathbb{R}^d))|\sqrt{t}f\in L^{\infty}(0,T;L^\infty(\mathbb{R}^d))\mathrm{\, and\, }\lim_{t\to 0}\|\sqrt{t}f\|_{L^{\infty}(\mathbb{R}^d)}=0\}$, and $\mathcal{F}_{T}=\{f|t^{1/4}f\in L^{\infty}(0,T;L^{2d,\infty}(\mathbb{R}^d))\mathrm{\, and\, } \lim_{t\to 0}\|t^{1/4}f\|_{L^{2d,\infty}(\mathbb{R}^d)}=0\}$. 
Then we can construct a mild solution $u=e^{t\Delta}u_{0}+A(f)-B(u,u)\in \mathcal{E}_{T_{0}}\subset\mathcal{F}_{T_{0}}$ for some positive $T_0$, where $A(f)=\int_{0}^{t}e^{(t-s)\Delta}\mathbb{P}fds$, $B(u,v)=\int_{0}^{t}e^{(t-s)\Delta}\mathbb{P}\nabla\cdot(u\otimes v)ds$. As proven in \cite[Chpater 15]{lemarie2002recent}, $A(f), B(u,u)\in C([0,T_0];L_{\sigma}^{d,\infty})$. Moreover, $A(f), B(u,u),e^{t\Delta}u_{0}\in C([0,T_{0}];\widetilde{L}_{\sigma}^{d,\infty})$ since $A(f), B(u,u)\in \mathcal{E}_{T_{0}}$ and $u_{0}\in X^{d,\infty}_{\sigma}$. This completes the proof. $\hfill\square$
\section*{Acknowledgement}
The author would like to express his gratitude to Professor Yasunori Maekawa, who instructs the author in the uniqueness problem and gives lots of valuable suggestions. 

\end{document}